\documentclass[reqno,12pt]{amsart}

\usepackage{amsmath,amsthm,amssymb,amsfonts,amscd,comment}
\usepackage{xypic}

\usepackage{graphicx}
\input xy
\xyoption{all}

\theoremstyle{plain} 
\newtheorem{thm}{Theorem}[section]
\newtheorem*{thm*}{Theorem}
\newtheorem{cor}[thm]{Corollary}
\newtheorem{lem}[thm]{Lemma}
\newtheorem{prop}[thm]{Proposition}

\newtheorem*{conj*}{Conjecture}

\theoremstyle{definition}
\newtheorem{defn}[thm]{Definition}%[section]
\newtheorem{eg}[thm]{Example}

\newtheorem{question}[thm]{Question}

\theoremstyle{remark}
\newtheorem{rem}[thm]{Remark}

\newtheorem*{pf}{Proof}

\numberwithin{equation}{section}
\def\ZZ{{\mathbb Z}}

\def\RR{{\mathbb R}}
\def\CC{{\mathbb C}}

\def\HH{{\mathbb H}}

\def\PP{{\mathbb P}}

\def\A{{\mathcal A}}

\def\C{{\mathcal C}}
\def\D{{\mathcal D}}

\def\N{{\mathcal N}}
\def\O{{\mathcal O}}
\def\P{{\mathcal P}}

\def\T{{\mathcal T}}

\newcommand{\Aut}{{\rm Aut}}
\newcommand{\Hom}{{\rm Hom}}

\newcommand{\re}{{\rm Re}\hspace{0.5mm}}
\newcommand{\im}{{\rm Im}\hspace{0.5mm}}

\def \mf#1#2#3#4{
\xymatrix{
{#1}\  \ar@<0.4ex>[r]^{{#2}} & \ {#4}
\ar@<0.4ex>[l]^{{#3}}
}
}

\def \mfs#1#2#3#4{\!
\xymatrix@C=1,5em{{#1} \! \ar@<0.2ex>[r]^{{#2}} & \! {#4}
\ar@<0.2ex>[l]^{{#3}}
}
\!}

\def \mfl#1#2#3#4{
\xymatrix@C=2.6em{{#1}\  \ar@<0.4ex>[r]^{{#2}} &\  {#4}
\ar@<0.2ex>[l]^{{#3}}
}
}

\def \mfss#1#2#3#4{\!
\xymatrix@C=1.5em{{#1} \ar@<0.3ex>[r]^{{#2}} & {#4}
\ar@<0.3ex>[l]^{{#3}}
}
\!}
\begin{document}
\title{Curvature of the space of stability conditions}
\date{\today}
\author{Kohei Kikuta}
\address{Department of Mathematics, Graduate School of Science, Osaka University, 
Toyonaka Osaka, 560-0043, Japan}
\email{k-kikuta@cr.math.sci.osaka-u.ac.jp}
\begin{abstract}
Motivated by the study of the autoequivalence group of triangulated categories via isometric actions on metric spaces, 
we consider curvature properties (CAT(0), Gromov hyperbolic) of the space of Bridgeland stability conditions with the canonical metric defined by Bridgeland. 
We then prove that the metric is neither CAT(0) nor hyperbolic, 
and the quotient metric by the natural $\CC$-action is not CAT(0) in case of the Kronecker quiver. 
Moreover, we also show the hyperbolicity of pseudo-Anosov functors defined by Dimitrov--Haiden--Katzarkov--Kontsevich, which yields the lower-bound of entropy by the translation length. 
Finally, pseudo-Anosov functors in case of curves have been completely classified. 
\end{abstract}
\maketitle
%\markboth{KOHEI KIKUTA}{TITLE}
%\tableofcontents
%%%%%%%%%%%%%%%%%%%%%%%%%%%%%%%%%%%%%%%%%%
\section{Introduction}
Let $\Gamma$ be a group. 
To study various properties of $\Gamma$, it is useful to construct isometric actions of $\Gamma$ on metric spaces with non-positive curvature property: CAT(0), hyperbolic in the sense of Gromov (cf. \cite[Part I\hspace{-0.1mm}I\hspace{-0.1mm}I.$\Gamma$]{BriH}). 
As an example, we consider the mapping class group ${\rm MCG}(M)$ of an oriented closed $2$-dimensional smooth manifold $M$. 
The group ${\rm MCG}(M)$ acts by isometries on the Teichm\"uller space with the Weil--Petersson metric, and the curve complex with a natural metric. 
The Weil--Petersson metric is CAT(0) and the curve complex is hyperbolic. 
%where the Weil--Petersson metric is CAT(0) and the curve complex is hyperbolic. 

%
We then focus on the autoequivalence group $\Aut(\T)$ of a triangulated category $\T$ (e.g. the derived category $\D^b(X)$ of coherent sheaves on a smooth projective variety $X$). 
Bridgeland introduced the notion of stability conditions on $\T$, and proved the set ${\rm Stab}(\T)$ of stability conditions admits a structure of complex manifolds (\cite{Bri1}). 
A metric $d_B$ on ${\rm Stab}(\T)$ was also introduced by Bridgeland, and then the group $\Aut(\T)$ naturally acts by isometries on $({\rm Stab}(\T),d_B)$. 
The group ${\rm MCG}(M)$ acts canonically on the (derived) Fukaya category of $M$, 
and the space of Bridgeland stability conditions is closely related to the K\"ahler moduli space or the moduli of complex structures (cf. \cite[Section 1]{Bri1}). 
Therefore, motivated by the homological mirror symmetry, it is natural to consider the analogue of various results on the isometric actions of ${\rm MCG}(M)$ for the isometric action of $\Aut(\T)$ on $({\rm Stab}(\T),d_B)$ (see also \cite{Smi}). 
Recently, Fan--Kanazawa--Yau actually defined an analogue of the Weil--Petersson metric on (the $\CC$-quotient of some subspace of) ${\rm Stab}(\D^b(X))$ for K3 surfaces $X$ (\cite{FKY}).

Although a motivation for the non-positivity of metrics is to study actions of autoequivalence groups in this paper, non-positive curvature metrics, especially CAT(0) metrics, are usually treated in the context of the contractibility of the space of stability conditions since any CAT(0) space is contractible (cf. \cite{BB,Smi}). 
Actually for K3 surface $X$, 
Allcock (\cite{All}) and Bridgeland (\cite{Bri1}) conjecture 
the CAT(0) property of 
the pull-back metric via the natural covering (different from Bridgeland's metric) on the $\CC$-quotient of some subspace of ${\rm Stab}(\D^b(X))$, see also \cite{Smi}. 

%----------
In this paper, we study the non-positive curvature property of $d_B$. 
Let $\T$ be a triangulated category of finite type over a field, and ${\rm Stab}^\dagger(\T)$ a connected component of the space ${\rm Stab}(\T)$ of stability conditions.  
The following is the first result: 
\begin{thm*}[Theorem \ref{d_B-thm}]
We have the followings: 
\begin{enumerate}
\item
The metric space $({\rm Stab}^{\dagger}(\T), d_B)$ is not CAT(0). 
\item
The metric space $({\rm Stab}^{\dagger}(\T), d_B)$ is not hyperbolic. 
\end{enumerate}
\end{thm*}
The key of the proof is to consider an orbit of the natural $\CC$-action on ${\rm Stab}^\dagger(\T)$. 
We can also show the similar statement (Theorem \ref{K3-red-d_B-thm}) for the space of reduced numerical stability conditions on the derived category of K3 surfaces of Picard rank one in the same way. 

Some examples tell us that the $\CC$-quotient space ${\rm Stab}^\dagger(\T)/\CC$ often admits a good geometric property. 
The Weil--Petersson metric due to Fan--Kanazawa--Yau (\cite{FKY}) is one of examples. 
We also consider the $\CC$-quotient metric $\bar{d_B}$ in the case of the derived category $\D^b({\rm mod}\hspace{0.5mm}\CC K_l)$ of the $l$-Kronecker quiver $K_l$, and prove the following: 
\begin{thm*}[Theorem \ref{quotient-d_B-Kronecker}]
The metric space $({\rm Stab}^{\dagger}(\D^b({\rm mod}\hspace{0.5mm}\CC K_l))/\CC, \bar{d_B})$ is not CAT(0). 
\end{thm*}
We also show the similar statement (Theorem \ref{Ginzburg-non-CAT(0)}) for the derived category of the Ginzburg differential graded $\CC$-algebra associated to $K_l$. 

The pseudo-Anosov functor introduced by Dimitrov--Haiden--Katzarkov--Kontsevich (\cite{DHKK}), is a categorical analogue of the pseudo-Anosov class in ${\rm MCG}(M)$. 
The following result is the categorical analogue of classical facts: the hyperbolicity of pseudo-Anosov classes with respect to the Teichm\"uller metric on the Teichm\"uller space, and the equality between the translation length and the stretch-factor for pseudo-Anosov classes. 
\begin{thm*}[Theorem \ref{pA-hyperbolic}]
Let 
%$\T$ be a triangulated category of finite type over a field, and 
$F\in\Aut(\T)$ be a pseudo-Anosov functor. 
Then, 
\begin{enumerate}
\item
$F$ is hyperbolic isometry with respect to the action on $({\rm Stab}^{\dagger}(\T)/\CC, \bar{d_B})$. 
\item
We have $\bar{d_B}(F)=\log\rho_F$, 
\end{enumerate}
where $\bar{d_B}(F)$ is the translation length of $F$ and $\rho_F$ is the stretch factor of $F$. 
\end{thm*}
As a corollary, we have the lower-bound of the entropy of $F$ by the translation length of $F$ (Corollary \ref{entropy-ge-translation}). 
In particular, for the derived category of curves, the pseudo-Anosov functors are completely classified and the entropy is equal to the translation length. 

%----------
The contents of this paper are as follows. 
In Section 2, we prepare some notations and define Bridgeland stability conditions, the metric $d_B$ on ${\rm Stab}(\T)$, natural group actions on ${\rm Stab}(\T)$ and the notions of CAT(0) and hyperbolic spaces. 
In Section 3, we study non-positive curvature properties of the metric $d_B$ and its $\CC$-quotient metric $\bar{d_B}$, and prove Theorem \ref{d_B-thm}, Theorem \ref{K3-red-d_B-thm}, Theorem \ref{quotient-d_B-Kronecker} and Theorem \ref{Ginzburg-non-CAT(0)}. 
In Section 4, based on \cite{DHKK}, pseudo-Anosov functors are introduced. 
We then prove Theorem \ref{pA-hyperbolic} and Corollary \ref{entropy-ge-translation}, and consider the case of curves. 
Some natural questions about the metric $d_B$ are listed in the final section. 

%-----------------------------------------------------------------------------------------------------------

\bigskip
\noindent
{\it Acknowledgment}\\
\indent
This work was partially done at The University of Sheffield which the author visited in the program ``Overseas Challenge Program for Young Researchers'' supported by JSPS, from September 2018 to February 2019. 
The author would like to express his sincere gratitude to the host Professor Tom Bridgeland for his kind hospitality and valuable discussions in Sheffield. 
The author would like to thank Genki Ouchi for telling him Proposition \ref{K3-pA} and its proof and Atsushi Takahashi for useful comments. 
The author is also supported by JSPS KAKENHI Grant Number JP17J00227. 
%%%%%%%%%%%%%%%%%%%%%%%%%%%%%%%%%%%%%%%%%%
\section{Preliminaries}
%-----------------------------------------------------------------------------------------------------------
\subsection{Notations}
Throughout this paper, $\T$ is a triangulated category of finite type over a field, and $\Aut(\T)$ is the group of autoequivalences of $\T$.  
We denote the Grothendieck group of $\T$ by $K(\T)$. 
%We denote the Grothendieck group (resp. numerical Grothendieck group) of $\T$ by $K(\T)$ (resp. $\N(\T)$). 
The numerical Grothendieck group $\N(\T)$ is defined as the quotient  of $K(\T)$ by the left radical of the Euler form. 

%-----------------------------------------------------------------------------------------------------------
\subsection{Stability condition}
In this subsection, we give the definition of Bridgeland stability conditions and basic properties of the canonical metric on the space of stability conditions. 

Fix a finitely generated free abelian group $\Lambda$, a surjective group homomorphism $v: K(\T)\twoheadrightarrow\Lambda$ and a norm $||\cdot||$ on $\Lambda\otimes_\ZZ \RR$. 
We assume moreover the existence of a group homomorphism $\alpha: \Aut(\T)\to\Aut_{\ZZ}(\Lambda)$, such that the following diagram of abelian groups is commutative for all $F\in\Aut(\T)$: 
\[
  \xymatrix{
    K(\T) \ar[r]^{K(F)} \ar[d]_{v} & K(\T) \ar[d]_{v} \\
    \Lambda \ar[r]^{\alpha(F)} & \Lambda 
  }
\]
Note that such $\alpha$, if exists, is uniquely determined. 

\begin{defn}[{\cite[Definition 5.1]{Bri1}}]%def of stab cond
A {\it stability condition} $\sigma=(Z,\P)$ on $\T$ (with respect to $(\Lambda, v)$) consists of a group homomorphism
$Z: \Lambda\to\CC$ called {\it central charge} and a family $\P=\{\P(\phi)\}_{\phi\in\RR}$ of full additive subcategory of $\T$ called {\it slicing}, such that
\begin{enumerate}
\item
For $0\neq E\in\P(\phi)$, we have $Z(v(E))=m(E)\exp(i\pi\phi)$ for some $m(E)\in\RR_{>0}$. 
\item
For all $\phi\in\RR$, we have $\P(\phi+1)=\P(\phi)[1]$.
\item
For $\phi_1>\phi_2$ and $E_i\in\P(\phi_i)$, we have $\Hom(E_1,E_2)= 0$.
\item
For each $0\neq E\in\T$, there is a collection of exact triangles called {\it Harder--Narasimhan filtration} of $E$: 
\begin{equation}\label{HN}
\begin{xy}
(0,5) *{0=E_0}="0", (20,5)*{E_{1}}="1", (30,5)*{\dots}, (40,5)*{E_{p-1}}="k-1", (60,5)*{E_p=E}="k",
(10,-5)*{A_{1}}="n1", (30,-5)*{\dots}, (50,-5)*{A_{p}}="nk",
\ar "0"; "1"
\ar "1"; "n1"
\ar@{.>} "n1";"0"
\ar "k-1"; "k" 
\ar "k"; "nk"
\ar@{.>} "nk";"k-1"
\end{xy}
\end{equation}
with $A_i\in\P(\phi_i)$ and $\phi_1>\phi_2>\cdots>\phi_p$. 
\item
(support property)
%We fix a norm $||\cdot||$ on $\Lambda\otimes_\ZZ \RR$. 
There exists a constant $C>0$ such that for all $0\neq E\in\P(\phi)$, 
we have 
\begin{equation}
||v(E)||<C|Z(v(E))|.
\end{equation}
\end{enumerate}
\end{defn}
For any interval $I\subset\RR$, define $\P(I)$ to be the extension-closed subcategory of $\T$ generated by the subcategories $\P(\phi)$ for $\phi\in I$. 
Then $\P((0,1])$ is the heart of a bounded t-structure on $\T$, hence an abelian category.  
The full subcategory $P(\phi)\subset\T$ is also shown to be abelian. 
A non-zero object $E\in\P(\phi)$ is called {\it $\sigma$-semistable} of {\it phase} $\phi_\sigma(E):=\phi$, and especially a simple object in $\P(\phi)$ is called {\it $\sigma$-stable}. 
Taking the Harder--Narasimhan filtration (\ref{HN}) of $E$, we define $\phi^+_\sigma(E):=\phi_\sigma(A_1)$ and $\phi^-_\sigma(E):=\phi_\sigma(A_p)$. 
The object $A_i$ is called {\it $\sigma$-semistable factor} of $E$. 
Define ${\rm Stab}_\Lambda(\T)$ to be the set of stability conditions on $\T$ with respect to $(\Lambda,v)$, 
especially, 
when $K(\T)$ (resp. $\N(\T)$) is finitely generated and free, 
${\rm Stab}_K(\T)$ (resp. ${\rm Stab}_\N(\T)$) to be the set of stability conditions on $\T$ with respect to $(K(\T),{\rm id})$ (resp. the natural projection $K(\T)\to\N(\T)$). 
An element in ${\rm Stab}_\N(\T)$ is called {\it numerical stability condition}. 

In this paper, we assume that the space ${\rm Stab}_\Lambda(\T)$ is not empty for some $(\Lambda,v)$. 
We will abuse notation and write $Z(E)$ instead of $Z(v(E))$. 

%---
\begin{comment}%Toda's remark
The above ‘good’ conditions are called ‘numerical property’ and ‘support property’ in literatures. Although the above properties are important in considering the space Stab(D), we omit the detail since we focus on the construction of one specific stability condition. 
\end{comment}

%-----------------------------------------------------------------------------------------------------------
We prepare some terminologies on the stability on the heart of a $t$-structure on $\T$. 
\begin{defn}%def of stability function
Let $\A$ be the heart of a bounded $t$-structure on $\T$. 
%and $K(\A)(\simeq\hspace{-0.2mm}K(\T))$ be its Grothendieck group. 
A {\it stability function} on $\A$ is a group homomorphism $Z: \Lambda\to\CC$ such that for all $0\neq E\in\A\subset\T$, the complex number $Z(v(E))$ lies in the semiclosed upper half plane $\HH_{-}:=\{re^{i\pi\phi}\in\CC~|~r\in\RR_{>0},\phi\in(0,1]\}\subset\CC$. 
\end{defn}
%We note that the heart of a bounded $t$-structure is an abelian category. 
Given a stability function $Z: \Lambda\to\CC$ on $\A$, the {\it phase} of an object $0\neq E\in\A$ is defined to be $\phi(E):=\frac{1}{\pi}{\rm arg}Z(E)\in(0,1]$. 
An object $0\neq E\in\A$ is {\it $Z$-semistable} (resp. {\it $Z$-stable}) if for all subobjects $0\neq A\subset E$, we have $\phi(A)\le\phi(E)$ (resp. $\phi(A)<\phi(E)$). 
We say that a stability function $Z$ satisfies {\it the Harder--Narasimhan property} if 
each object $0\neq E\in\A$ admits a filtration (called Harder--Narasimhan filtration of $E$) 
$0=E_0\subset E_1\subset E_2\subset\cdots\subset E_m=E$ such that $E_i/E_{i-1}$ is $Z$-semistable for $i=1,\cdots,m$ with $\phi(E_1/E_0)>\phi(E_2/E_1)>\cdots>\phi(E_m/E_{m-1})$, 
and {\it the support property} if there exists a constant $C>0$ such that for all $Z$-semistable objects $E\in\A$, 
we have $||v(E)||<C|Z(v(E))|$. 

The following proposition shows the relationship between stability conditions and stability functions on the heart of a bounded $t$-structure. 
\begin{prop}[{\cite[Proposition 5.3]{Bri1}}]
To give a stability condition on $\T$ is equivalent to giving the heart $\A$ of a bounded t-structure on $\T$, and a stability function $Z$ on $\A$ with the Harder--Narasimhan property and the support property. 
\end{prop}
For the proof, we construct the slicing $\P$, from the pair $(Z,\A)$, by 
\[
\P(\phi):=\{E\in\A~|~E\text{ is }Z\text{-semistable with }\phi(E)=\phi\}\text{ for }\phi\in(0,1],
\]
and extend for all $\phi\in\RR$ by $\P(\phi+1):=\P(\phi)[1]$. 
Conversely, for a stability condition $\sigma=(Z,\P)$, the heart $\A$ is given by $\A:=\P_\sigma((0,1])$. 
We also denote stability conditions by $(Z,\A)$. 

%----------
The following notion is important to analyze the space of stability conditions. 
\begin{defn}%def of mass
Let $E\in\T$ be a non-zero object of $\T$ and $\sigma\in {\rm Stab}_\Lambda(\T)$ be a stability condition on $\T$. 
The {\it mass} $m_{\sigma}(E)\in\RR_{>0}$ of $E$ is defined by
\begin{equation}
m_{\sigma}(E):=\displaystyle\sum_{i=1}^p|Z_\sigma(A_i)|, 
\end{equation}
where $A_1,\cdots,A_p$ are $\sigma$-semistable factors of $E$. 
\end{defn}

The following generalized metric (i.e. with values in $[0,\infty]$) defined by Bridgeland plays a central role in this paper.   
\begin{defn}[{\cite[Proposition 8.1]{Bri1}}]%def of d_B
The generalized metric $d_B$ on ${\rm Stab}_\Lambda(\T)$ is defined by 
\begin{equation}
d_B(\sigma, \tau):=\sup_{E\neq0}\left\{|\phi^+_\sigma(E)-\phi^+_\tau(E)|, |\phi^-_\sigma(E)-\phi^-_\tau(E)|, \middle|\log\frac{m_\sigma(E)}{m_\tau(E)}\middle|\right\}\in[0,\infty]. 
%d_B(\sigma, \tau):=\sup_{E\neq0}\left\{\max\left\{|\phi^+_\sigma(E)-\phi^+_\tau(E)|, |\phi^-_\sigma(E)-\phi^-_\tau(E)|, \middle|\log\frac{m_\sigma(E)}{m_\tau(E)}\middle|\right\}\right\}\in[0,\infty]. 
\end{equation}
\end{defn}
This generalized metric induces the topology on ${\rm Stab}_\Lambda(\T)$. 
Then the generalized metric $d_B$ takes a finite value on each connected component ${\rm Stab}_\Lambda^\dagger(\T)$ of ${\rm Stab}_\Lambda(\T)$, 
thus $({\rm Stab}_\Lambda^\dagger(\T),d_B)$ is a metric space in the strict sense. 
\begin{thm}[{\cite[Theorem 7.1]{Bri1}}]%local homeo
The map
\begin{equation}
{\rm Stab}_\Lambda(\T)\to\Hom_\ZZ(\Lambda, \CC);~\sigma=(Z,\P)\mapsto Z
\end{equation}
is a local homeomorphism, where $\Hom_\ZZ(\Lambda, \CC)$ is equipped with the natural linear topology. 
\end{thm}
Therefore the space ${\rm Stab}_\Lambda(\T)$ (and each connected component ${\rm Stab}_\Lambda^\dagger(\T)$) naturally admits a structure of finite dimensional complex manifolds and is especially locally compact. 
\begin{thm}[{\cite[Thorem 3.6]{Woo}}]
The metric space $({\rm Stab}_\Lambda^{\dagger}(\T), d_B)$ is complete. 
Moreover, the limit point $\sigma_\infty=(Z_\infty, \P_\infty)$ of a Cauchy sequence $\{\sigma_n=(Z_n,\P_n)\}_n\subset {\rm Stab}_\Lambda^\dagger(\T)$ is described by
\begin{eqnarray*}
&&Z_\infty=\lim_{n\to\infty}Z_n,\\
&&\P_\infty(\phi)=\{0\}\cup\langle0\neq E\in\T~|~\lim_{n\to\infty}\phi^+_{\sigma_n}(E)=\phi,~\lim_{n\to\infty}\phi^-_{\sigma_n}(E)=\phi\rangle. 
\end{eqnarray*}
\end{thm}
There are two group-actions on ${\rm Stab}_\Lambda(\T)$. 
The first is a left $\Aut(\T)$-action defined by 
\begin{equation}
F.\sigma:=(Z_\sigma(\alpha(F^{-1})v(-)), \{F(\P_\sigma(\phi))\})~\text{ for }\sigma\in {\rm Stab}_\Lambda(\T),~F\in\Aut(\T). 
\end{equation}
Let $\varphi: \widetilde{{\rm GL}}_+(2,\RR)\to{\rm GL}_+(2,\RR)$ be the universal cover of ${\rm GL}_+(2,\RR)$. 
We recall that $\widetilde{{\rm GL}}_+(2,\RR)$ is isomorphic to the group of pairs $(M,f)$ where $f:\RR\to\RR$ is an increasing map with $f(\phi+1)=f(\phi)+1$, and $M\in{\rm GL}_+(2,\RR)$ such that the induced maps on $(\RR^2\backslash\{0\})/\RR_{>0}=S^1=\RR/2\ZZ$ coincide. 
The second is a right $\widetilde{{\rm GL}}_+(2,\RR)$-action defined by
\begin{equation}\label{univ-cov-gl}
\sigma.g:=(M^{-1}\circ Z_\sigma, \{\P_\sigma(f(\phi))\})~\text{ for }\sigma\in {\rm Stab}_\Lambda(\T),~g=(M,f)\in\widetilde{{\rm GL}}_+(2,\RR). 
\end{equation}
These two actions commute, and the $\widetilde{{\rm GL}}_+(2,\RR)$-action is free and continuous. 
Restricting the $\widetilde{{\rm GL}}_+(2,\RR)$-action to the subgroup $\CC\subset\widetilde{{\rm GL}}_+(2,\RR)$, 
we also have a right $\CC$-action defined by
\begin{equation}
\sigma.\lambda:=(\exp(-i\pi\lambda)\cdot Z_\sigma, \{\P_\sigma(\phi+\re\lambda)\})\text{ for }\lambda\in\CC. 
\end{equation}
\begin{lem}
The $\Aut(\T)$-action and $\CC$-action are isometry with respect to $d_B$. 
\end{lem}
\begin{comment}%---comment out in ver1
Associated to a Riemannian metric on a $C^\infty$-manifold, the length metric is canonically defined by the infimum of the (Riemannian) length of piecewise differentiable paths. 
The following proposition has already mentioned in \cite[Section 2]{Woo}. 
\begin{prop}
There is no Riemannian metric on the complex manifold ${\rm Stab}_\Lambda^\dagger(\T)$ whose length metric is isometric to $d_B$. 
\end{prop}
\begin{pf}
\qed
\end{pf}
\end{comment}
%-----------------------------------------------------------------------------------------------------------
\subsection{Metric space}
We collect basic definitions in metric geometry. 
We refer to the text \cite{BriH} and it is highly recommended for details.  
For a metric space $(X,d)$, the open ball with center $x\in X$ and radius $r>0$ is denoted by $B(x,r)$. 
Let $(\RR^2,d_E)$ be the Euclidean plane. 
\begin{defn}%def of geodesic, geodesic space, uniquely geodesic
Let $(X,d)$ be a metric space. 
\begin{enumerate}
\item
Let $x,y\in X$ be two distinct points. 
A continuous map $\gamma: [0,1]\to X$ satisfying $\gamma(0)=x$ and $\gamma(1)=y$ is a {\it geodesic} from $x$ to $y$ if $d(\gamma(t),\gamma(t'))=d(x,y)\cdot|t-t'|$ for all $t,t'\in[0,1]$. 
\item
A metric space $(X,d)$ is {\it geodesic} if, for all two distinct points $x,y\in X$, there exists a geodesic from $x$ to $y$. 
\item%def of uniquely geodesic
A metric space $(X,d)$ is {\it uniquely geodesic} if, for all two distinct points $x,y\in X$, there uniquely exists a geodesic from $x$ to $y$. 
\end{enumerate}
We denote a geodesic from $x$ to $y$ (or its image in $X$) by $[x,y]$. 
This standard notation is less accurate if $(X,d)$ is not uniquely geodesic, but we use it for simplicity.  
\end{defn}
\begin{defn}%comparison triangle, comparison point
Let $(X,d)$ be a metric space. 
\begin{enumerate}
\item
A {\it geodesic triangle} $\Delta([x,y],[y,z],[z,x]):=[x,y]\cup[y,z]\cup[z,x]\subset X$ is the union of geodesics $[x,y],[y,z],[z,x]$. 
\item
A geodesic triangle $\Delta(\bar{x},\bar{y},\bar{z})$ in $(\RR^2,d_E)$ is a {\it comparison triangle} for a geodesic triangle $\Delta([x,y],[y,z],[z,x])$ in $X$ if
\begin{equation}
d(x,y)=d_E(\bar{x},\bar{y}),~d(y,z)=d_E(\bar{y},\bar{z}),~d(z,x)=d_E(\bar{z},\bar{x}).  
\end{equation}
\item
Fix a comparison triangle $\Delta(\bar{x},\bar{y},\bar{z})$ for a geodesic triangle $\Delta([x,y],[y,z],[z,x])$ in $X$. 
A point $\bar{p}\in[\bar{x},\bar{y}]$ is a {\it comparison point} for $p\in[x,y]$ if $d(x,p)=d_E(\bar{x},\bar{p})$. 
%which is uniquely determined if we fix a comparison triangle. 
Comparison points on $[\bar{y},\bar{z}], [\bar{z},\bar{x}]$ are defined in the same way. 
\end{enumerate}
\end{defn}
For a geodesic triangle in a metric space $(X,d)$, a comparison triangle uniquely exists up to isometry (see \cite[Lemma I.2.14]{BriH}). 
We then define the CAT(0) property of metric spaces. 
\begin{defn}%def of CAT(0)
A metric space $(X,d)$ is {\it CAT(0)} if 
\begin{enumerate}
\item
$(X,d)$ is geodesic. 
\item
For each geodesic triangle $\Delta([x,y],[y,z],[z,x])$ in $X$, we fix a comparison triangle $\Delta(\bar{x},\bar{y},\bar{z})$. 
Then, for all $p,q\in\Delta([x,y],[y,z],[z,x])$, we have 
\begin{equation}
d(p,q)\le d_E(\bar{p},\bar{q}), 
\end{equation}
where $\bar{p},\bar{q}\in\Delta(\bar{x},\bar{y},\bar{z})$ are comparison points for $p,q\in\Delta([x,y],[y,z],[z,x])$. 
\end{enumerate}
\end{defn}
We call a metric space $(X,d)$ {\it locally CAT(0)} if, for each point $x\in X$, there exists an open neighborhood of $x$ such that the induced metric is CAT(0). 
Note that CAT(0) spaces are uniquely geodesic, locally CAT(0) and contractible (see \cite[Proposition II.1.4]{BriH}). 

We then define the hyperbolicity (in the sense of Gromov) of metric spaces. 
\begin{defn}%def of \delta-hyperbolic
Fix $\delta\in\RR_{>0}$. 
A metric space $(X, d)$ is {\it $\delta$-hyperbolic} if 
\begin{enumerate}
\item
$(X,d)$ is geodesic. 
\item
All geodesic triangles $\Delta([x,y],[y,z],[z,x])$ in $X$ satisfy
\begin{equation}
[x,y]\subset N_\delta([y,z]\cup[z,x]),~[y,z]\subset N_\delta([z,x]\cup[x,y]),~[z,x]\subset N_\delta([x,y]\cup[y,z]), 
\end{equation}
where, for a subset $Y\subset X$, $N_\delta(Y)\subset X$ is the $\delta$-neighborhood of $Y$. 
%\delta-slim(not \delta-thin)
\end{enumerate}
\end{defn}
\begin{defn}%def of hyperbolic
A metric space $(X, d)$ is {\it hyperbolic} if $(X, d)$ is $\delta$-hyperbolic for some $\delta\in\RR_{>0}$. 
\end{defn}
The following lemma is easy to show, but we note that this is used in the proof of main theorems. 
\begin{lem}\label{subset}
Let $(X,d)$ be a metric space and $Y$ a subset of X. 
\begin{enumerate}
\item
If $(Y,d|_Y)$ is geodesic and not $CAT(0)$, then $(X,d)$ is not $CAT(0)$.  
\item
If $(Y,d|_Y)$ is geodesic and not hyperbolic, then $(X,d)$ is not  hyperbolic. 
\end{enumerate}
\end{lem}
%%%%%%%%%%%%%%%%%%%%%%%%%%%%%%%%%%%%%%%%%%
\section{Curvature properties}
In this section we study the non-positive curvature property (i.e. CAT(0), hyperbolic) of the metric $d_B$ and the $\CC$-quotient metric $\bar{d_B}$, and prove Theorem \ref{d_B-thm} and Theorem \ref{quotient-d_B-Kronecker}. 
\subsection{The metric $d_B$}
Let ${\rm Stab}_\Lambda^{\dagger}(\T)$ be a connected component of ${\rm Stab}_\Lambda(\T)$, 
and fix a stability condition $\sigma\in{\rm Stab}_\Lambda^{\dagger}(\T)$.  
We now study metric properties of $\CC$-orbits. 
\begin{lem}\label{C-d_B}
For each $\lambda\in\CC$, we have
\begin{equation}
d_B(\sigma,\sigma.\lambda)=\max\{|\re\lambda|, \pi|\im\lambda|\}. 
\end{equation}
\end{lem}
\begin{pf}
It is clear since the $\CC$-action preserves semistable objects and Harder--Narasimhan filtrations. 
\qed
\end{pf}
\begin{lem}[{\cite[Section 2]{Woo}}]\label{C-orbit-geod}
The metric space $(\sigma.\CC, d_B)$ is geodesic. 
\end{lem}
\begin{pf}
For any $\lambda,\lambda'\in\CC$, the straight line $\gamma(t):=\sigma.((1-t)\lambda+t\lambda')$ for $t\in[0,1]$ is a geodesic from $\sigma.\lambda$ to $\sigma.\lambda'$. 
\qed
\end{pf}
\begin{prop}\label{d_B-curvature}
We have the followings: 
\begin{enumerate}
\item
$(\sigma.\CC, d_B)$ is not locally CAT(0) (especially not CAT(0)). 
\item
$(\sigma.\CC, d_B)$ is not hyperbolic. 
\end{enumerate}
\end{prop}
\begin{pf}
The straight line as in the proof of Lemma \ref{C-orbit-geod} from $x$ to $y$ $(x,y\in\sigma.\CC)$ is denoted by $[x,y]$. 

\noindent
(i) It is enough to show that $(B(\sigma, r')\cap\sigma.\CC, d_B)$ is not uniquely geodesic for all $r'\in\RR_{>0}$ since all CAT(0) spaces are uniquely geodesic. 
For $r\in(0,r')\subset\RR_{>0}$, we set $\lambda:=r, \lambda_1:=\frac{r}{2}(1+\frac{i}{2\pi})$. 
%\lambda_2:=\frac{r}{2}(1-\frac{i}{2\pi})$. 
Then we have 
\begin{equation}
d_B(\sigma,\sigma.\lambda)=d_B(\sigma,\sigma.\lambda_1)+d_B(\sigma.\lambda_1,\sigma.\lambda). 
\end{equation}
It is easy to see that geodesics $[\sigma,\sigma.\lambda], [\sigma,\sigma.\lambda_1], [\sigma.\lambda_1,\sigma.\lambda]$ are contained in $B(\sigma, r')\cap\sigma.\CC$ and $\sigma.\lambda_1\not\in[\sigma,\sigma.\lambda]$. 
Thus the (reparametrization of) path $[\sigma,\sigma.\lambda_1]\ast[\sigma.\lambda_1,\sigma.\lambda]$ is a geodesic from $\sigma$ to $\sigma.\lambda$, which is not equal to $[\sigma,\sigma.\lambda]$. 
Therefore $(B(\sigma, r')\cap\sigma.\CC, d_B)$ is not uniquely geodesic. 
%since there are two different geodesics $[\sigma,\sigma.\lambda],[\sigma,\sigma.\lambda_1]\cup[\sigma.\lambda_1,\sigma.\lambda]$ from $\sigma$ to $\sigma.\lambda$. 

%reparameterization is also geodesic which differs from the straight line. 
%Therefore, $(\sigma.\CC, d_B|_{\sigma.\CC})$ is not uniquely geodesic. 

%---
\noindent
(ii) Fix $\delta>0$. 
We set $x_1:=\sigma, x_2:=\sigma.4\delta,x_3:=\sigma.(i\frac{4\delta}{\pi}),x_4:=\sigma.(2\delta+i\frac{2\delta}{\pi})$ and let $\Delta_\delta:=\Delta(\gamma_{12},\gamma_{23},\gamma_{13})$ be a geodesic triangle on $\sigma.\CC$, 
where $\gamma_{ij}:=[x_i,x_j]$. 
It is clear that $x_4\in\gamma_{23}([0,1])$
Then, for all $t\in[0,1]$, we have
\begin{equation}
d_B(\gamma_{12}(t),x_4)
=\max\{|4t-2|\delta,2\delta\}
\ge 2\delta>\delta
\end{equation}
and similarly we get $d_B(\gamma_{13}(t),x_4)>\delta$. 
We thus have $\gamma_{23}\not\subset N_\delta(\gamma_{12}\cup\gamma_{13})$ i.e. $\Delta_\delta$ is not $\delta$-slim. 
Therefore $(\sigma.\CC, d_B)$ is not $\delta$-hyperbolic for all $\delta>0$, hence not hyperbolic. 
\qed
\end{pf}
The following is one of the main theorems in this section: 
\begin{thm}\label{d_B-thm}
We have the followings: 
\begin{enumerate}
\item
$({\rm Stab}_\Lambda^{\dagger}(\T), d_B)$ is not CAT(0). 
\item
$({\rm Stab}_\Lambda^{\dagger}(\T), d_B)$ is not hyperbolic. 
\end{enumerate}
\end{thm}
\begin{pf}
We fix a stability condition $\sigma\in {\rm Stab}_\Lambda^{\dagger}(\T)$, then the claim follows from Lemma \ref{C-orbit-geod}, Proposition \ref{d_B-curvature} and Lemma \ref{subset}. 
\qed
\end{pf}
%-----------------------------------------------------------------------------------------------------------
\subsubsection{K3 surfaces}\label{K3}
Let $S$ be a complex algebraic K3 surface of Picard rank $1$ and $\D^b(S)$ the bounded derived category of coherent sheaves on $S$. 
The space of stability conditions in this case is important in the context of mirror symmetry (cf. \cite{Bri1,BB}). 
Let ${\rm Stab}_\N^\dagger(\D^b(S))\subset {\rm Stab}_\N(\D^b(S))$ be the connected component containing the set of geometric stability conditions, and ${\rm Stab}^\dagger_{red}(\D^b(S))$ the subset of ${\rm Stab}_\N^\dagger(\D^b(S))$ consisting of reduced stability conditions (see \cite{Bri2, BB} for definitions). 
The space ${\rm Stab}_\N^\dagger(\D^b(S))$ is contractible by the celebrated theorem due to Bayer--Bridgeland (\cite[Theorem 1.3]{BB}). 
The following is clear since the metric space $({\rm Stab}^\dagger_{red}(\D^b(S)), d_B)$ is preserved by the $\CC$-action. 
\begin{thm}\label{K3-red-d_B-thm}
The metric space $({\rm Stab}^\dagger_{red}(\D^b(S)), d_B)$ is neither CAT(0) nor hyperbolic. 
\end{thm}
%-----------------------------------------------------------------------------------------------------------
\subsection{The $\CC$-quotient metic $\bar{d_B}$}
In this subsection, we shall consider the quotient, by the $\CC$-action, of some connected component of the space of stability conditions. 

Let $\bar{d_B}$ be the quotient metric of $d_B$ by the $\CC$-action, that is, $\bar{d_B}$ is defined by
\begin{equation}
\bar{d_B}(\bar{\sigma},\bar{\tau}):=\inf_{\lambda\in\CC}d_B(\sigma,\tau.\lambda)
\end{equation}

It is easy to see that every $\CC$-orbit is a path-connected closed subset in ${\rm Stab}_\Lambda(\T)$ by Lemma \ref{C-d_B} and Lemma \ref{C-orbit-geod}. 
Thus the quotient metric $\bar{d_B}$ is well-defined on ${\rm Stab}^\dagger_\Lambda(\T)/\CC$. 

%\begin{eg}
%Let $Q_{A_1}$ be the $A_1$-quiver and $D^b({\rm mod}\hspace{0.5mm}\CC Q_{A_1})$ the bounded derived category of finitely generated $\CC Q_{A_1}(\simeq\CC)$-modules. 
%Then ${\rm Stab}_K(\D^b({\rm mod}\hspace{0.5mm}\CC Q_{A_1}))/\CC$ is an one point set. 
%\end{eg}
Let $\HH:=\{x+iy\in\CC~|~x,y\in\RR,~y>0\}$ be the upper half plane and $d_P$ the ${\rm Poincar\acute{e}}$ metric on $\HH$ associated to the Riemannian metric 
$g_P=\frac{dx^2+dy^2}{4y^2}$. 
\begin{eg}
Let $C$ be a smooth projective curve over $\CC$ of genus $g(C)\ge1$ and $\D^b(C)$ the bounded derived category of coherent sheaves on $C$. 
Then ${\rm Stab}_\N(\D^b(C))$ is connected and, by \cite[Proposition 4.1]{Woo}, the quotient metric space $({\rm Stab}_\N(\D^b(C))/\CC, \bar{d_B})$ is a proper $CAT(0)$ hyperbolic space since there is an isometry $({\rm Stab}_\N(\D^b(C))/\CC, \bar{d_B})\simeq(\HH, d_P)$. 
\end{eg}
%-----------------------------------------------------------------------------------------------------------
\subsubsection{Kronecker quiver}
Let $Q$ be a connected finite acyclic quiver and $\D^b({\rm mod}\hspace{0.5mm}\CC Q)$ the bounded derived category of finitely generated $\CC Q$-modules. 
The simple module corresponding to a vertex $i$ of $Q$ is denoted by $S_i\in{\rm mod}\hspace{0.5mm}\CC Q$. 
A group homomorphism $Z: K(\D^b({\rm mod}\hspace{0.5mm}\CC Q))\to\CC$ such that $Z(S_i)\in\HH_-$ for each vertex $i$ is a stability function on $\D^b({\rm mod}\hspace{0.5mm}\CC Q)$, and obviously satisfies the support property. 
Then the pair $\sigma_Z:=(Z,{\rm mod}\hspace{0.5mm}\CC Q)$ becomes a stability condition on $\D^b({\rm mod}\hspace{0.5mm}\CC Q)$ by \cite[Proposition 2.4]{Bri1}. 
We define the distinguished connected component ${\rm Stab}_K^\dagger(\D^b({\rm mod}\hspace{0.5mm}\CC Q))$ of ${\rm Stab}_K(\D^b({\rm mod}\hspace{0.5mm}\CC Q))$ to be the connected component containing the connected subset
\begin{equation}
{\rm Stab}(Q):=\{\sigma_Z=(Z,{\rm mod}\hspace{0.5mm}\CC Q)\}\subset {\rm Stab}_K(\D^b({\rm mod}\hspace{0.5mm}\CC Q)). 
\end{equation}

%---Kronecker quiver

Let $K_l$ be the $l$-Kronecker quiver for $l\in\ZZ_{>0}$ as follows: 
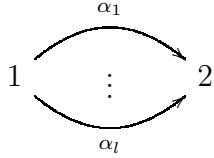
\[
\xymatrix{
1 \ar@(ur,ul)[rr]^{\alpha_1} \ar@{}[rr]|{\vdots} \ar@(dr,dl)[rr]_{\alpha_l} &  & 2
}
 \]
The following is one of the main theorems in this section: 
\begin{thm}\label{quotient-d_B-Kronecker}
The metric space $({\rm Stab}_K^{\dagger}(\D^b({\rm mod}\hspace{0.5mm}\CC K_l))/\CC, \bar{d_B})$ is not CAT(0). 
\end{thm}
To prove this theorem, we shall construct a geodesic, but non-CAT(0) subspace in ${\rm Stab}_K^\dagger(\D^b({\rm mod}\hspace{0.5mm}\CC K_l))/\CC$ for applying  Lemma \ref{subset}. 
\begin{defn}
We define the metric $d'_B$ on $\RR^4$ by 
\begin{equation}
d'_B(x,x'):=\max_{i=1,\cdots,4}\{|x_i-x'_i|\}
\end{equation}
for $x=(x_i)_{i=1,\cdots,4}, x'=(x'_i)_{i=1,\cdots,4}\in\RR^4$, 
and a (right) $\CC$-action on $\RR^4$ by
\begin{equation}
x.\lambda:=(x_1+\re\lambda, x_2+\im\lambda, x_3+\re\lambda, x_4+\im\lambda)
\end{equation}
for $x=(x_i)_{i=1,\cdots,4}\in\RR^4$ and $\lambda\in\CC$. 
\end{defn}
For any $x,x'\in\RR^4$, the straight line $\gamma(t):=(1-t)x+tx'$ for $t\in[0,1]$ is a geodesic, 
thus $(\RR^4,d'_B)$ is a geodesic space. 
The above $\CC$-action is obviously isometric with respect to $d'_B$. 
The quotient metric $\bar{d'_B}$ on $\RR^4/\CC$ is well-defined since each $\CC$-orbit is closed. 
\begin{lem}
The quotient metric $\bar{d'_B}$ on $\RR^4/\CC$ is described as follows: 
\begin{equation}
\bar{d'_B}(\bar{x},\bar{x'})=\max\left\{\frac{|(x'_1-x_1)-(x'_3-x_3)|}{2},~\frac{|(x'_2-x_2)-(x'_4-x_4)|}{2}\right\}
\end{equation}
for all $\bar{x},\bar{x'}\in\RR^4/\CC$. 
Moreover, $(\RR^4/\CC, \bar{d'_B})$ is geodesic. 
\end{lem}
\begin{pf}
The first claim follows from
\begin{eqnarray*}
d'_B(x.\lambda,x')
&=&\inf_{\lambda\in\CC}\max\{|\re\lambda-(x'_1-x_1)|, |\im\lambda-(x'_2-x_2)|, |\re\lambda-(x'_3-x_3)|, |\im\lambda-(x'_4-x_4)|\}\\
&\ge&\max\left\{\frac{|(x'_1-x_1)-(x'_3-x_3)|}{2},~\frac{|(x'_2-x_2)-(x'_4-x_4)|}{2}\right\}\\
&=&d'_B\left(x.\left(\frac{(x'_1-x_1)+(x'_3-x_3)}{2}+i\frac{(x'_2-x_2)+(x'_4-x_4)}{2}\right),x'\right). 
\end{eqnarray*}
The quotient of straight lines in $\RR^4$ are also geodesics in $\RR^4/\CC$. 
\qed
\end{pf}
%def of the chamber $\C$
The subset $\C\subset\RR^4$ is defined by 
\begin{equation}
\C:=\{x\in\RR^4~|~0<x_1<x_3\le1\}, 
\end{equation}
and set $\overline{\C}:=(\C.\CC)/\CC\subset\RR^4/\CC$. 
\begin{lem}\label{bar-C-geod}
The metric space $(\overline{\C}, \bar{d'_B})$ is geodesic. 
\begin{pf}
By the definition of $\C$ it is easy to see that, for $x,x'\in\C.\CC$, the straight line from $x$ to $x'$ is contained in $\C.\CC$. 
The claim follows from that the quotient of straight lines are geodesics in $\overline{\C}$. 
\qed
\end{pf}
\end{lem}
\begin{prop}\label{bar-C-nonCAT(0)}
The metric space $(\overline{\C}, \bar{d'_B})$ is not locally CAT(0) (especially not CAT(0)). 
\end{prop}
\begin{pf}
The quotient of straight line in $\overline{\C}$ from $\bar{x}$ to $\bar{y}$ $(x,y\in\C)$ is denoted by $[x,y]$. 
For any $r'\in\RR_{>0}$, we fix $r\in\RR_{>0}$ such that $0<r<\max\{r', \frac{1}{4}\}$, 
and set
\begin{equation}
p_1:=\overline{(r,0,2r,0)},~p_2:=\overline{\left(r,0,3r,\frac{r}{2}\right)},~p_3:=\overline{(r,0,4r,0)}\in\overline{\C}. 
\end{equation}
Then it is enough to show that $(B(p_1, r')\cap\overline{\C}, \bar{d'_B})$ is not uniquely geodesic. 
We have 
\begin{equation}
\bar{d_B'}(p_1,p_3)=\bar{d_B'}(p_1,p_2)+\bar{d_B'}(p_2,p_3), 
\end{equation}
It is easy to see that geodesics $[p_1,p_3], [p_1,p_2],[p_2,p_3]$ are contained in $B(p_1, r')\cap\overline{\C}$ and $p_2\not\in[p_1,p_3]$. 
Thus the (reparametrization of) path $[p_1,p_2]\ast[p_2,p_3]$ is a geodesic from $p_1$ to $p_3$, which is not equal to $[p_1,p_3]$. 
Therefore $(B(p_1, r')\cap\overline{\C}, \bar{d'_B})$ is not uniquely geodesic. 
\qed
\end{pf}
%-----------------------------------------------------------------------------------------------------------

\vspace{3mm}
\noindent{\it Proof of Theorem \ref{quotient-d_B-Kronecker}}
\vspace{2mm}

\noindent 
%-----------------------------------------------------------------------------------------------------------
We define a map $q: \C.\CC\to {\rm Stab}(K_l).\CC\subset {\rm Stab}_K^{\dagger}(\D^b({\rm mod}\hspace{0.5mm}\CC K_l))$ by
\begin{equation}
q(x.\lambda):=(Z,{\rm mod}\hspace{0.5mm}\CC K_l).\left(\re\lambda+i\frac{\im\lambda}{\pi}\right), 
\end{equation}
where the stability function $Z$ on ${\rm mod}\hspace{0.5mm}\CC K_l$ is defined by
\begin{equation}
Z(S_1):=\exp(x_2+i\pi x_1), ~Z(S_2):=\exp(x_4+i\pi x_3)\in\HH_-. 
\end{equation}
This map is well-defined by direct computations. 
Then we will show that $q$ is an isometric embedding. 
%----------q: isom emb
For each $\sigma=(Z,{\rm mod}\hspace{0.5mm}\CC K_l)\in q(\C)$, we have $0<\phi_\sigma(S_1)<\phi_\sigma(S_2)\le1$. 
Thus the Harder--Narasimhan filtration of each object $0\neq E\in {\rm mod}\hspace{0.5mm}\CC K_l$ is the form of 
\begin{equation}
\begin{xy}
(0,5) *{0}="u0", (20,5)*{S_2^{\oplus k_2}}="u1", (40,5)*{E}="u2",
(10,-5)*{S_2^{\oplus k_2}}="d1", (30,-5)*{S_1^{\oplus k_1}}="d2",
\ar "u0"; "u1"
\ar "u1"; "d1"
\ar@{.>} "d1";"u0"
\ar "u1"; "u2" 
\ar "u2"; "d2"
\ar@{.>} "d2";"u1"
\end{xy}, 
\end{equation}
that is, for all $\sigma\in q(\C)$, Harder--Narasimhan filtrations of each nonzero object in $\D^b({\rm mod}\hspace{0.5mm}\CC K_l)$ are the same. 
Then for all $\sigma=q(x),\tau=q(y)\in q(\C)$ , we have 
\begin{eqnarray*}
d_B(q(x),q(y))
&=&\sup_{E\neq0}\left\{|\phi^+_\sigma(E)-\phi^+_\tau(E)|, |\phi^-_\sigma(E)-\phi^-_\tau(E)|, \middle|\log\frac{m_\sigma(E)}{m_\tau(E)}\middle|\right\}\\
&=&\max_{i=1,2}\left\{|\phi_\sigma(S_i)-\phi_\tau(S_i)|, \middle|\log\frac{|Z_\sigma(S_i)|}{|Z_\tau(S_i)|}\middle|\right\}\\
&=&\max_{j=1,\cdots,4}\{|x_j-y_j|\}\\
&=&d_B'(x,y), 
\end{eqnarray*}
where we use an elementary inequality 
$\left|\log\frac{r_1+r_2}{r_1'+r_2'}\right|\le\max\left\{\left|\log\frac{r_1}{r_1'}\right|,\left|\log\frac{r_2}{r_2'}\right|\right\}$ for all $r_1,r_1',r_2,r_2'\in\RR_{>0}$
 to estimate the mass term. 
Therefore $q$ is an isometric embedding on $\C$. 
Since the $\CC$-action preserves semistable objects and Harder--Narasimhan filtrations, similar computations imply that the map $q$ is an isometric embedding (on $\C.\CC$). 

%----------
The induced map
\begin{equation}
\bar{q}: \overline{\C}\to {\rm Stab}(K_l).\CC/\CC\subset {\rm Stab}_K^{\dagger}(\D^b({\rm mod}\hspace{0.5mm}\CC K_l))/\CC
\end{equation}
is also well-defined by the definition of $q$. 
Since $q$ is an isometric embedding, it is easy to check that $\bar{q}$ is also isometric embedding. 
Therefore the proof is completed by Lemma \ref{bar-C-geod}, Proposition \ref{bar-C-nonCAT(0)} and Lemma \ref{subset}. 
\qed
\vspace{3mm}

%-----------------------------------------------------------------------------------------------------------
The complex manifold ${\rm Stab}_K^{\dagger}(\D^b({\rm mod}\hspace{0.5mm}\CC K_l))/\CC$ is biholomorphic to $\CC$ for $l=1,2$ (cf. \cite{Qiu, Oka}) and $\HH$ for $l\ge3$ (cf. \cite[Theorem 1.5, Corollary 3.3]{DK}). 
On the other hand, as a corollary of Theorem \ref{quotient-d_B-Kronecker}, these biholomorphisms are not isometries with respect to their standard metrics: 
\begin{cor}
We have the followings: 
\begin{enumerate}
\item
$({\rm Stab}^{\dagger}(\D^b({\rm mod}\hspace{0.5mm}\CC K_l))/\CC, \bar{d_B})$ is not isometric to $(\CC, d_E)$ for $l=1,2$. 
\item
$({\rm Stab}^{\dagger}(\D^b({\rm mod}\hspace{0.5mm}\CC K_l))/\CC, \bar{d_B})$ is not isometric to $(\HH, d_P)$ for $l\ge3$.
\end{enumerate}
\end{cor}
\begin{pf}
It is obvious since the Euclidean metric $d_E$ and the ${\rm Poincar\acute{e}}$ metric $d_P$ is CAT(0). 
\qed
\end{pf}
Let $\Gamma_N K_l$ be the Ginzburg $N$-Calabi--Yau differential graded (dg) $\CC$-algebra associated to $K_l$ for $N\ge3$ (cf. \cite[Section 4.2]{Gin} or \cite[Section 6.2]{Kel}), and $\D_{fd}(\Gamma_N K_l)$ the derived category of finite dimensional dg $\Gamma_N K_l$-modules. 
There canonically exists a bounded $t$-structure on $\D_{fd}(\Gamma_N K_l)$ whose heart is ${\rm mod}\hspace{0.5mm}\CC K_l$ (cf. \cite[Section 2]{Ami}). 
Thus, by the same arguments as in the case of $\D^b({\rm mod}\hspace{0.5mm}\CC K_l)$, we can construct the isometric embedding
\begin{equation}
\overline{\C}\to {\rm Stab}^{\dagger}_K(\D_{fd}(\Gamma_N K_l))/\CC, 
\end{equation}
which yields the following: 
\begin{thm}\label{Ginzburg-non-CAT(0)}
The metric space $({\rm Stab}^{\dagger}_K(\D_{fd}(\Gamma_N K_l))/\CC, \bar{d_B})$ is not CAT(0). 
\end{thm}
%-----------------------------------------------------------------------------------------------------------

%%%%%%%%%%%%%%%%%%%%%%%%%%%%%%%%%%%%%%%%%%
\section{Hyperbolicity of pseudo-Anosov functors}
In this section, we consider the isometric action of $\Aut(\T)$ on ${\rm Stab}_\Lambda(\T)$, especially the hyperbolicity of a certain isometry called ``pseudo-Anosov functor''. 
%-----------------------------------------------------------------------------------------------------------
\subsection{Hyperbolic isometry}
We shall recall the notion of hyperbolic isometry. 
Let $(X,d)$ be a metric space and $\Gamma$ a group acting by isometries on $(X,d)$. 
\begin{defn}%translation length
We define the {\it translation length} $d(g)$ of $g\in\Gamma$ with respect to the isometric action on $(X,d)$ by
\begin{equation}
d(g):=\displaystyle\inf_{x\in X}d(x,g.x)\in\RR_{\ge0}. 
\end{equation}
\end{defn}
\begin{defn}
Let $(X, d)$ be a metric space. 
An isometry $g\in\Gamma$ of $(X, d)$ is {\it hyperbolic} 
if there exists $x\in X$ such that $d(g)=d(x,g.x)$ and $d(g)>0$. 
%If the translation length $d(g)$ is attained and $d(g)>0$, then $f$ is called {\it hyperbolic}. 
\end{defn}
By the triangle inequality, we have 
\begin{equation}
d(x,g.x)\ge\displaystyle\limsup_{n\to\infty}\frac{1}{n}d(x,g^n.x)
\end{equation}
for all $x\in X$. 
Using the triangle inequality again, the right hand side of the above inequality is independent of the choice of $x$, hence \begin{equation}
d(g)\ge\displaystyle\limsup_{n\to\infty}\frac{1}{n}d(x,g^n.x). 
\end{equation}
Since the sequence $\{d(x,g^n.x)\}_n$ satisfies the subadditivity: 
$d(x,g^{m_1+m_2}.x)\le d(x,g^{m_1}.x)+d(x,g^{m_2}.x)$ for all $m_1,m_2\in\ZZ_{>0}$, 
the limit $\lim_{n\to\infty}\frac{1}{n}d(x,g^n.x)$ exists by the classical fact, thus we have 
\begin{equation}
d(g)\ge\displaystyle\lim_{n\to\infty}\frac{1}{n}d(x,g^n.x). 
\end{equation}
%---
\begin{comment}
It is well-known that
the inequality 
\begin{equation}
d(g)\ge\displaystyle\lim_{n\to\infty}\frac{1}{n}d(x,g^n.x)
\end{equation}
always holds for all $x\in X$ and each isometry $g\in\Gamma$ by the triangle inequality. 
\end{comment}
%---
This actually becomes an equality when $(X,d)$ is CAT(0) and $g$ is a hyperbolic isometry (cf. \cite[Theorem II.6.8]{BriH}). 
% and each point $x\in X$ on any axis of $g$: Hence any point in $X$. 
%-----------------------------------------------------------------------------------------------------------
\subsection{Mass-growth and entropy}
We collect basics on dynamical invariants on $\Aut(\T)$. 
\begin{defn}[{\cite[Section 4]{DHKK} and \cite[Theorem 3.5(1)]{Ike}}]%def of mass-growth
Let $G\in\T$ be a split-generator, $F\in\Aut(\T)$ an autoequivalence of $\T$ and 
$\sigma\in {\rm Stab}_\Lambda(\T)$ a stability condition on $\T$. 
The {\it mass-growth} $h_{\sigma}(F)\in[-\infty,+\infty]$ with respect to $F$ is defined by
\begin{equation}
h_{\sigma}(F):=\displaystyle\limsup_{n\rightarrow\infty}\frac{1}{n}\log m_{\sigma}(F^n G), 
\end{equation}
\end{defn}
In the definition of the mass-growth, the limit doesn't depend on the choice of split-generators (see \cite[Theorem 3.5(1)]{Ike}). 
%---
\begin{comment}
\begin{prop}
For every autoequivalence $F\in\Aut^\dagger(\T)$, we have 
\begin{equation}
d_B(F)\ge\displaystyle\lim_{n\to\infty}\frac{1}{n}d_B(\sigma, F^n\sigma)\ge h_\sigma(F)
\end{equation}
\end{prop}
\begin{pf}
\qed
\end{pf}
\end{comment}
%---
\begin{defn}[{\cite[Definition 2.5]{DHKK}}]%def of entropy
Let $G\in\T$ be a split-generator and $F\in\Aut(\T)$. 
The {\it entropy} $h(F)$ of $F$ is defined as follows:  
\begin{equation}
h(F):=\displaystyle\lim_{n\rightarrow\infty}\frac{1}{n}\log \delta(G,F^{n}G)\in\RR_{\geq0}, 
\end{equation}
where
the complexity $\delta(G, F^nG)$ is the minimum $p\in\ZZ_{>0}$ such that there exists a diagram of exact triangles of the following form
\begin{equation}
\begin{xy}
(0,5) *{0}="0", (20,5)*{M_{1}}="1", (30,5)*{\dots}, (40,5)*{M_{p-1}}="k-1", (60,5)*{F^nG\oplus M}="k",
(10,-5)*{G[n_{1}]}="n1", (30,-5)*{\dots}, (50,-5)*{G[n_{p}]}="nk",
\ar "0"; "1"
\ar "1"; "n1"
\ar@{.>} "n1";"0"
\ar "k-1"; "k" 
\ar "k"; "nk"
\ar@{.>} "nk";"k-1"
\end{xy}
\end{equation}
\begin{comment}
\begin{equation}
\delta(G,F^{n}G)
:=
\min\left\{
\displaystyle p\in\ZZ_{>0}
~\middle|~
\begin{xy}
(0,5) *{0}="0", (20,5)*{M_{1}}="1", (30,5)*{\dots}, (40,5)*{M_{p-1}}="k-1", (60,5)*{F^nG\oplus M}="k",
(10,-5)*{G[n_{1}]}="n1", (30,-5)*{\dots}, (50,-5)*{G[n_{p}]}="nk",
\ar "0"; "1"
\ar "1"; "n1"
\ar@{.>} "n1";"0"
\ar "k-1"; "k" 
\ar "k"; "nk"
\ar@{.>} "nk";"k-1"
\end{xy}
\right\}
\end{equation}
\end{comment}
\end{defn}
In the definition of the entropy, the limit exists and doesn't depend on the choice of split-generators (see \cite[Lemma 2.6]{DHKK}). 
\begin{thm}[{\cite[Thorem 3.5(2)]{Ike}}]\label{mass-growth-le-entropy}
Let $F\in\Aut(\T)$ be an autoequivalence of $\T$ and $\sigma\in {\rm Stab}_\Lambda(\T)$ a stability condition on $\T$. 
Then we have
\begin{equation}
h_{\sigma}(F)\le h(F)<\infty.
\end{equation}
\end{thm}
%-----------------------------------------------------------------------------------------------------------
\subsection{Pseudo-Anosov functor}
%and $\Aut^\dagger(\T)$ the subgroup of $\Aut(\T)$ which preserving ${\rm Stab}^\dagger_\Lambda(\T)$. 
The pseudo-Anosov functor is introduced due to Dimitrov--Haiden--Katzarkov--Kontsevich, motivated by a categorical analogue of pseudo-Anosov mapping class acting on the (derived) Fukaya category of Riemann surfaces. 

We fix a connected component ${\rm Stab}_\Lambda^{\dagger}(\T)$ of ${\rm Stab}_\Lambda(\T)$. 
%Let ${\rm Stab}_\Lambda^{\dagger}(\T)$ be a connected component of ${\rm Stab}_\Lambda(\T)$. 
\begin{defn}[{\cite[Definition 4.1]{DHKK}}]\label{pA}%def of pA
An autoequivalence $F\in\Aut(\T)$ is {\it pseudo-Anosov} if 
there exist a stability condition $\sigma_F\in {\rm Stab}_\Lambda^\dagger(\T)$ and an element $g_F\in\widetilde{{\rm GL}}_+(2,\RR)$ such that 
\begin{enumerate}
\item
$\varphi(g_F)=
\begin{pmatrix}
\frac{1}{r}&0\\
0&r
\end{pmatrix}
\text{ or }
\begin{pmatrix}
r&0\\
0&\frac{1}{r}
\end{pmatrix}
\in {\rm GL}_+(2,\RR)$ for some $|r|>1$. 
\item
$F.\sigma_F=\sigma_F.g_F$
\end{enumerate}
We here recall $\varphi: \widetilde{{\rm GL}}_+(2,\RR)\to{\rm GL}_+(2,\RR)$ is the universal cover of ${\rm GL}_+(2,\RR)$, and call $\rho_F:=|r|>1$ the stretch-factor of $F$. 
\end{defn}
Pseudo-Anosov functors preserve the connected component ${\rm Stab}_\Lambda^\dagger(\T)$ since the $\widetilde{{\rm GL}}_+(2,\RR)$-action is continuous. 
\begin{eg}[{\cite[Section 4]{DHKK}}]
Let $\D^b({\rm mod}\hspace{0.5mm}\CC K_l)$ be the derived category of the $l$-Kronecker quiver $K_l$ for $l\ge3$. 
Then the Serre functor of $\D^b({\rm mod}\hspace{0.5mm}\CC K_l)$ is pseudo-Anosov. 
\end{eg}
The mass-growth with respect to the pseudo-Anosov is described as follows. 
\begin{prop}\label{pA-mass-growth}
For a pseudo-Anosov functor $F\in\Aut(\T)$, we have
\begin{equation}
h_{\sigma_F}(F)=\log\rho_F>0. 
\end{equation}
\end{prop}
\begin{pf}
We fix a split-generator $G\in\T$ and let $\{A_i\}_{i=1,\cdots,p}$ be the $\sigma_F$-semistable factors of $G$ with $\phi_{\sigma_F}(A_{i-1})>\phi_{\sigma_F}(A_i)$. 
The $\sigma_F$-semistable factors of $FG$ are $\{FA_i\}_{i=1,\cdots,p}$ with $\phi_{\sigma_F}(FA_{i-1})>\phi_{\sigma_F}(FA_i)$
since pseudo-Anosov functors preserve $\sigma_F$-semistable objects and Harder--Narasimhan filtrations, and we have $m_{\sigma_F}(FA_i)=m_{F^{-1}.\sigma_F}(A_i)=m_{\sigma_F.g_F^{-1}}(A_i)$. 
Then the claim follows from a direct computation of the mass-growth. 
\qed
\end{pf}
As a corollary, we see that the stretch-factor is uniquely determined since the mass-growth is invariant on the same connected component (\cite[Proposition 3.10]{Ike}). 

%---
We prove the categorical analogue of classical facts: the hyperbolicity of pseudo-Anosov classes with respect to the Teichm\"uller metric, and the equality between the translation length and the stretch-factor for pseudo-Anosov classes. 
This is the main result in this section: 
\begin{thm}\label{pA-hyperbolic}
Let $F\in\Aut(\T)$ be a pseudo-Anosov functor. 
Then, 
\begin{enumerate}
\item
$F$ is hyperbolic isometry with respect to the action on $({\rm Stab}_\Lambda^{\dagger}(\T)/\CC, \bar{d_B})$. 
\item
We have $\bar{d_B}(F)=\displaystyle\lim_{n\to\infty}\frac{1}{n}\bar{d_B}(\bar{\sigma_F}, F^n\bar{\sigma_F})=\log\rho_F$. 
\end{enumerate}
\end{thm}
\begin{pf}
We set $\sigma:=\sigma_F$ and $g:=g_F$ for simplicity. 
Note that the $\widetilde{{\rm GL}}_+(2,\RR)$-action preserves semistable objects and Harder--Narasimhan filtrations. 
For all $\lambda\in\CC$, we have 
\begin{eqnarray*}
d_B(\sigma,F\sigma\lambda)
&=&d_B(\sigma,\sigma(g\circ\lambda))
=\sup_{E\neq0}\left\{|\phi^\pm_{\sigma(g\circ\lambda)}(E)-\phi^\pm_\sigma(E)|, \middle|\log\frac{m_{\sigma(g\circ\lambda)}(E)}{m_\sigma(E)}\middle|\right\}\\
&=&\sup_{E:~\sigma{\text -semistable}}\left\{|\phi_{\sigma(g\circ\lambda)}(E)-\phi_\sigma(E)|, \middle|\log\frac{|Z_{\sigma(g\circ\lambda)}(E)|}{|Z_\sigma(E)|}\middle|\right\}\\
&=&\sup_{E:~\sigma{\text -semistable}}\left\{|f_{g\circ\lambda}(\phi_\sigma(E))-\phi_\sigma(E)|, \middle|\log \left(\middle|M_{g\circ\lambda}\frac{Z_\sigma(E)}{|Z_\sigma(E)|}\middle|\right)\middle|\right\}\\
&=&\max\left\{\sup_{E:~\sigma{\text -semistable}}|f_{g\circ\lambda}(\phi_\sigma(E))-\phi_\sigma(E)|, \sup_{E:~\sigma{\text -semistable}}\middle|\log \left(\middle|M_{g\circ\lambda}\frac{Z_\sigma(E)}{|Z_\sigma(E)|}\middle|\right)\middle|\right\}\\
&\le&\max\left\{\sup_{t\in\RR}|f_{g\circ\lambda}(t)-t|, \sup_{v\in\RR^2, ||v||_E=1}|\log \left(\middle|M_{g\circ\lambda}v\middle|\right)|\right\}\\
&=&\max\left\{||f_{g\circ\lambda}-{\rm id}_\RR||_\RR, \log||M_{g\circ\lambda}||, \log||M_{g\circ\lambda}^{-1}||\right\},
\end{eqnarray*}
where $||\cdot||_\RR$ is the sup-norm, $||\cdot||$ is the operator norm induced by the Euclidean norm $||\cdot||_E$ on $\RR^2$ and $g=(M_g,f_g)\in\widetilde{{\rm GL}}_+(2,\RR)$ (see (\ref{univ-cov-gl})). 
Moreover a direct computation implies that
\begin{eqnarray*}
& &\inf_{\lambda\in\CC}\max\{||f_{g\circ\lambda}-{\rm id}_\RR||,\log||M_{g\circ\lambda}||,\log||M_{g\circ\lambda}^{-1}||\}\\
&=&\max\{||f_g-{\rm id}_\RR||,\log||M_g||,\log||M_g^{-1}||\}\\
&=&\max\{\log||M_g||,\log||M_g^{-1}||\}
=\log\rho_F. 
\end{eqnarray*}
Therefore, we have 
\begin{equation}
\bar{d_B}(\bar{\sigma}, F\bar{\sigma})
\le\inf_{\lambda\in\CC}\max\{||f_{g\circ\lambda}-{\rm id}_\RR||,\log||M_{g\circ\lambda}||,\log||M_{g\circ\lambda}^{-1}||\}
=\log\rho_F. 
\end{equation}
See also the proof of \cite[Proposition 4.1]{Woo} for the above computations. 

Since for a $\sigma$-semistable object $E$, $F^k(E)$ and $F^{-k}(E)$ are also $\sigma$-semistable for $k\gg0$, we have 
%Considering $F^k(E)$ and $F^{-k}(E)$ for a $\sigma$-semistable object $E$ and $k\gg0$, we have
\begin{eqnarray}
\sup_{E:~\sigma{\text -semistable}}\left\{\log\middle|M_{g^n}\frac{Z_\sigma(E)}{|Z_\sigma(E)|}\middle|\right\}=\log\rho_F^n\\
\inf_{E:~\sigma{\text -semistable}}\left\{\log\middle|M_{g^n}\frac{Z_\sigma(E)}{|Z_\sigma(E)|}\middle|\right\}=\log\rho_F^{-n}. 
\end{eqnarray}
which implies
\begin{eqnarray*}
& &\inf_{\lambda\in\CC}\sup_{E:~\sigma{\text-semistable}}\left|\log \left(\middle|M_{g^n\circ\lambda}\frac{Z_\sigma(E)}{|Z_\sigma(E)|}\middle|\right)\right|\\
&=&\inf_{\lambda\in\CC}\sup_{E:~\sigma{\text-semistable}}\left|\log \left(\middle|M_{g^n}\frac{Z_\sigma(E)}{|Z_\sigma(E)|}\middle|\right)+\pi\im\lambda\right|\\
&=&\inf_{\lambda\in\CC}\max\{|\log\rho_F^n+\pi\im\lambda|,~|\log\rho_F^{-n}+\pi\im\lambda|\}
=\log\rho_F^n. 
\end{eqnarray*}
Therefore we have 
\begin{eqnarray*}
& &\bar{d_B}(\bar{\sigma}, F^n\bar{\sigma})\\
&=&\inf_{\lambda\in\CC}\max\left\{\sup_{E:~\sigma{\text -semistable}}|f_{g^n\circ\lambda}(\phi_\sigma(E))-\phi_\sigma(E)|, \sup_{E:~\sigma{\text -semistable}}\middle|\log \left(\middle|M_{g^n\circ\lambda}\frac{Z_\sigma(E)}{|Z_\sigma(E)|}\middle|\right)\middle|\right\}\\
&\ge&\inf_{\lambda\in\CC}\sup_{E:~\sigma{\text-semistable}}\left|\log \left(\middle|M_{g^n\circ\lambda}\frac{Z_\sigma(E)}{|Z_\sigma(E)|}\middle|\right)\right|=\log\rho_F^n, 
\end{eqnarray*}
which yields
\begin{equation}
\displaystyle\lim_{n\to\infty}\frac{1}{n}\bar{d_B}(\bar{\sigma}, F^n\bar{\sigma})
\ge \log\rho_F. 
\end{equation}
Summarizing the above argument, we have
\begin{eqnarray*}
\log\rho_F
\ge\bar{d_B}(\bar{\sigma}, F\bar{\sigma})
&\ge& \displaystyle\bar{d_B}(F)=\inf_{\bar{\sigma'}\in {\rm Stab}_\Lambda^{\dagger}(\T)/\CC}\bar{d_B}(\bar{\sigma'},F\bar{\sigma'})\\
&\ge& \displaystyle\lim_{n\to\infty}\frac{1}{n}\bar{d_B}(\bar{\sigma}, F^n\bar{\sigma})
\ge \log\rho_F, 
\end{eqnarray*}
which completes the proof. 
\qed
\end{pf}
\begin{rem}
If a stability condition $\sigma_F$ for a pseudo-Anosov functor $F\in\Aut(\T)$ has dense phase in $\RR$, 
then $d_B(\sigma_F,F\sigma_F)=\bar{d_B}({\bar\sigma_F},F\bar{\sigma_F})$ (See the proof of \cite[Proposition 4.1]{Woo}). 
\end{rem}
As a corollary, we have the lower-bound of entropy by the translation length: 
\begin{cor}\label{entropy-ge-translation}
For a pseudo-Anosov functor $F\in\Aut(\T)$, we have
\begin{equation}
h(F)\ge\bar{d_B}(F)>0. 
\end{equation}
\end{cor}
\begin{pf}
The claim follows from Theorem \ref{mass-growth-le-entropy}, Proposition \ref{pA-mass-growth} and Theorem \ref{pA-hyperbolic} (ii). 
\qed
\end{pf}
A stability condition is {\it algebraic} if 
its heart is a finite length abelian category with finitely many isomorphism classes of simple objects. 
A connected component of ${\rm Stab}_\Lambda(\T)$ containing an algebraic stability condition is called {\it algebraic}.  

\begin{cor}
For a pseudo-Anosov functor $F\in\Aut(\T)$ such that $\sigma_F$ is in an algebraic component ${\rm Stab}^{\dagger}_\Lambda(\T)$, 
we have
\begin{equation}
h(F)=\log\rho_F=\bar{d_B}(F). 
\end{equation}
\end{cor}
\begin{pf}
For $\sigma_F$ in the statement, we have $h(F)=h_{\sigma_F}(F)$ by \cite[Theorem 3.14]{Ike}.  
Therefore the claim follows from Proposition \ref{pA-mass-growth} and Theorem \ref{pA-hyperbolic} (ii). 
\qed
\end{pf}
%The following is an example of algebraic components: 
Let $Q$ be a connected finite acyclic quiver, and 
${\rm Stab}_K^\dagger(\D^b({\rm mod}\hspace{0.5mm}\CC Q))$ the distinguished component containing the set of algebraic stability conditions whose heart is ${\rm mod}\hspace{0.5mm}\CC Q$. 
The above corollary gives a generalization of DHKK's computation for the Serre functor of $\D^b({\rm mod}\hspace{0.5mm}\CC K_l)$ (cf. \cite[Section 4]{DHKK}). 

%Ikeda's result states that h(F) = h\UTF{2000}(F) whenever \UTF{2000} is in a component of the stability space containing an algebraic stability condition, i.e. one with a length heart with finitely many iso-classes of simple objects. Therefore Corollary 4.12 holds in greater generality.

%You could add this to Cor 4.11, by saying that h(F)\le \bar{d}_B(F) with equality when there is an algebraic stability condition in the component of $\sigma$, and then give the quiver case as an example.

\begin{comment}
Let $Q$ be a connected finite acyclic quiver, and 
recall that ${\rm Stab}_K^\dagger(\D^b(Q))$ is the distinguished component containing the set of stability conditions whose heart is ${\rm mod}\hspace{0.5mm}Q$. 
We can also show a generalization of DHKK's computation for the Serre functor of $\D^b(K_l)$ (cf. \cite[Section 4]{DHKK}) as a corollary: 
\begin{cor}
For a pseudo-Anosov functor $F\in\Aut(\D^b(Q))$ such that $\sigma_F\in {\rm Stab}_K^\dagger(\D^b(Q))$, we have
\begin{equation}
h(F)=\log\rho_F=\bar{d_B}(F). 
\end{equation}
\end{cor}
\begin{pf}
For $\sigma_F$ in the statement, we have $h(F)=h_{\sigma_F}(F)$ by \cite[Theorem 3.14]{Ike}.  
Therefore the claim follows from Proposition \ref{pA-mass-growth} and Theorem \ref{pA-hyperbolic} (ii). 
\qed
\end{pf}
\end{comment}

%-----------------------------------------------------------------------------------------------------------
\subsubsection{Curves}
We shall classify pseudo-Anosov functors for the derived category $\D^b(C)$ of coherent sheaves on a smooth projective curve $C$ over $\CC$. 
\begin{prop}\label{pA-not-ell}
Let $C$ be a smooth projective curve over $\CC$ with genus $g(C)\neq1$. 
Then no autoequivalence of $\D^b(C)$ is pseudo-Anosov. 
\end{prop}
\begin{pf}
The claim immediately follows from the fact that the entropy of each autoequivalence of curves with $g(C)\neq1$ is zero by \cite[Proposition 3.3]{Kik}. 
\qed
\end{pf}
Let $E$ be an elliptic curve with a closed point $x_0\in E$, and fix the basis $\{[\O_E],[\O_{x_0}]\}$ of $\N(\D^b(E))$. 
Then, via the isomorphism $\N(\D^b(E))\simeq\ZZ^2$, the automorphism $\N(F)$ on $\N(\D^b(E))$ induced by $F\in\Aut(\D^b(E))$ can be seen as an element in ${\rm SL}(2,\ZZ)$ since $\N(F)$ preserves the Euler form on $\N(\D^b(E))$. 

We can define a group homomorphism $Z_0: \N(\D^b(E))\to\CC$ by
\begin{equation}
Z_0(E):=-\mathrm{deg}(E)+i\cdot \mathrm{rk}(E), 
\end{equation}
which is a stability function on the standard heart ${\rm Coh}(E)$. 
%Let $\P_{\beta,H}$ is the slicing constructed by the $Z_{\beta,H}$-stability on ${\rm Coh}(C)$ which is the heart of the standard t-structure on $\D^b(C)$. 
Then the pair $\sigma_0:=(Z_0, {\rm Coh}(C))$ is a numerical stability condition on $\D^b(C)$. %(see \cite[Example 5.4]{Bri} for details)
The $\widetilde{{\rm GL}}_+(2,\RR)$-action on the space ${\rm Stab}_\N(\D^b(E))$ is free and transitive by \cite[Theorem 9.1]{Bri1}, 
thus ${\rm Stab}_\N(\D^b(E))$ is equal to $\sigma_0.\widetilde{{\rm GL}}_+(2,\RR)$. 
\begin{prop}\label{pA-ell}
An autoequivalence $F\in\Aut(\D^b(E))$ is psuedo-Anosov if and only if we have $|{\rm trace}(\N(F))|>2$. 
\end{prop}
\begin{pf}
For each autoequivalence $F'\in\Aut(\D^b(E))$ with $|{\rm trace}(\N(F'))|\le2$, we have $h(F')=0$ by \cite[Proposition 3.6]{Kik}, thus $F'$ is not pseudo-Anosov. 
Therefore it is enough to show that each autoequivalence $F\in\Aut(\D^b(E))$ with $|{\rm trace}(\N(F'))|>2$ is pseudo-Anosov. 
Note that $Z_0$ is given by the matrix
$
\begin{pmatrix}
0&-1\\
1&0
\end{pmatrix}
\in{\rm SL}(2,\ZZ)
$
after identifying $\N(\D^b(E))\simeq\ZZ^2$ using the chosen basis and $\CC\simeq\RR^2$ via real and imaginary parts. 
%Note that $Z_0$ can be identified with 
%via the fixed isomorphism $\N(\D^b(E))\simeq\ZZ^2$. 
By the equality ${\rm Stab}_\N(\D^b(E))=\sigma_0.\widetilde{{\rm GL}}_+(2,\RR)$, there exists $\tilde{g}\in\widetilde{{\rm GL}}_+(2,\RR)$ such that $F.\sigma_0=\sigma_0.\tilde{g}$,
which implies
$
\varphi(\tilde{g})=
\begin{pmatrix}
0&-1\\
1&0
\end{pmatrix}
\N(F)
\begin{pmatrix}
0&-1\\
1&0
\end{pmatrix}^{-1}
$. 
Thus we have $\varphi(\tilde{g})\in{\rm SL}(2,\ZZ)$ with $|{\rm trace}(\varphi(\tilde{g}))|>2$. 
%since det and trace are conjugation-invariant. 
The diagonalization of $\varphi(\tilde{g})$ yields 
$\varphi(\tilde{g})=hg_rh^{-1}$, where $h\in{\rm GL}_+(2,\RR)$ and $g_r:=
\begin{pmatrix}
\frac{1}{r}&0\\
0&r
\end{pmatrix}
\text{ or }
\begin{pmatrix}
r&0\\
0&\frac{1}{r}
\end{pmatrix}
$ 
with $|r|>1$. 
We can take lifts $\tilde{h}, \tilde{g_r}\in\widetilde{{\rm GL}}_+(2,\RR)$ of $h,g_r$ such that $\tilde{g}=\tilde{h}\tilde{g_r}\tilde{h}^{-1}$. 
Therefore $F$ is pseudo-Anosov since we have $F.(\sigma_0.\tilde{h})=(\sigma_0.\tilde{h}).\tilde{g_r}$
\qed
\end{pf}
Therefore we have completely classified pseudo-Anosov functors in the case of curves by Proposition \ref{pA-not-ell} and Proposition \ref{pA-ell}. 

Finally, the entropy of pseudo-Anosov functors of (elliptic) curves is described by the stretch-factor and the translation length with respect to $\bar{d_B}$, which is also categorical analogue of the classical fact: the topological entropy of pseudo-Anosov class is equal to the stretch-factor, hence the translation length. 
\begin{lem}
Let $E$ be an elliptic curve. 
For a pseudo-Anosov functor $F\in\Aut(\D^b(E))$, we have
\begin{equation}
h(F)=\log\rho_F=\bar{d_B}(F).
\end{equation}
\end{lem}
\begin{pf}
By \cite[Theorem 1.2]{Ike} and \cite[Theorem 3.11]{Kik}, we have $h(F)=h_{\sigma_F}(F)$. 
Therefore the claim follows from Proposition \ref{pA-mass-growth} and Theorem \ref{pA-hyperbolic} (ii). 
\qed
\end{pf}
%%%%%%%%%%%%%%%%%%%%%%%%%%%%%%%%%%%%%%%%%%
\section{Some questions}
In this section, we shall list some questions about the metric $d_B$. 
\begin{defn}%def of length metric
A metric space $(X,d)$ is a length space if we have 
\begin{equation}
d(x,y)=\inf_{\gamma}\ell(\gamma)~\text{ for }x,y\in X, 
\end{equation}
where $\gamma: [0,1]\to X$ is a continuous map such that $\gamma(0)=1$ and $\gamma(1)=y$
and $\ell(\gamma)$ is a length of $\gamma$ defined by  
\begin{equation}
\ell(\gamma):=\sup_{0=t_0\le t_1\le\cdots\le t_n=1}\sum_{i=0}^{n-1}d(\gamma(t_i),\gamma(t_{i+1})). 
\end{equation}
\end{defn}
We note that a geodesic space is a length space since the length of geodesics is equal to the distance. 
 The metric space $({\rm Stab}_\Lambda^{\dagger}(\T), d_B)$ is locally compact and complete, 
%It is easy to show that a proper metric space is complete and locally compact.  
thus the following holds by applying the Hopf--Rinow theorem (cf. \cite[Proposition I.3.7]{BriH}): 
\begin{prop}
If $({\rm Stab}_\Lambda^{\dagger}(\T), d_B)$ is a length space, then $({\rm Stab}_\Lambda^{\dagger}(\T), d_B)$ is a proper geodesic space. 
\end{prop}
Therefore the following question is natural: 
\begin{question}
The metric space $({\rm Stab}_\Lambda^{\dagger}(\T), d_B)$ is a length space?
\end{question}
We note that if $d_B$ is a length metric, then the quotient metric $\bar{d_B}$ is also length metric. 

%---
In the previous section, we have got the classification of pseudo-Anosov functors in the case of curves, 
while K3 surfaces doesn't admit pseudo-Anosov functors with respect to stability conditions in ${\rm Stab}_\N^\dagger(\D^b(S))$ (see \ref{K3} for notations). 
\begin{prop}[Ouchi]\label{K3-pA}
Let $S$ be a complex algebraic K3 surfaces. 
Then there is no pseudo-Anosov functor $F$ such that $\sigma_F\in{\rm Stab}_\N^\dagger(\D^b(S))$. 
\end{prop}
\begin{pf}
Let $(-,-)$ be the Mukai pairing on $H^*(S,\ZZ)$, and $v: \N(\D^b(S))\hookrightarrow H^*(S,\ZZ)$ the Mukai vector. 
Assume that $F\in\Aut(\D^b(S))$ is a pseudo-Anosov functor such that $\sigma_F=(Z_F,\P_F)\in{\rm Stab}_\N^\dagger(\D^b(S))$. 
Then there exists a vector $\Omega_F\in v(\N(\D^b(S)))\otimes_\ZZ\CC$ satisfying $Z_F(-)=(\Omega_F,v(-))$, and that $\re\Omega_F$ and $\im\Omega_F$ span a positive definite real plane in $v(\N(\D^b(S)))\otimes_\ZZ\RR$. 
Since the induced automorphism $F^H$ on $H^*(S,\CC)$ preserves $H^*(S,\ZZ)$ and the Mukai pairing and it follows that $Z_F(F^{-1}(-))=\varphi(g_F)^{-1}Z_F(-)$, 
we have 
\begin{equation}
F^H(\re\Omega_F)=r\re\Omega_F\text{  or  }\frac{1}{r}\re\Omega_F, 
\end{equation}
which contradicts to $(\re\Omega_F)^2>0$ by $|r|>1$. 
\qed
\end{pf}
Therefore we shall suggest the following question. 
%the following question is for the next step: 
\begin{question}
Let $S$ be a smooth projective surface over $\CC$. 
Then when does $\D^b(S)$ have a pseudo-Anosov functor?  
\end{question}

%---
Recently, Fan--Filip--Haiden--Katzarkov--Liu propose a generalization of pseudo-Anosov functor in \cite{FFHKL}. 
\begin{question}
The results of Section 4 can be generalized in the sense of \cite{FFHKL}?
%The recent preprint arXiv:1910.12350 by Fan, Filip, Haiden, Katzarkov and Liu proposes a more general de\UTF{2000}nition of pseudo-Anosov functor. It might be worth adding questions about whether the results of §4 generalise to their setting?
\end{question}
%%%%%%%%%%%%%%%%%%%%%%%%%%%%%%%%%%%%%%%%%%
%\appendix
%%%%%%%%%%%%%%%%%%%%%%%%%%%%%%%%%%%%%%%%%%

\end{document}